\documentclass[letter,final,oneside,notitlepage,11pt]{article}

\usepackage{amssymb}
\usepackage{amsmath}
\usepackage{amstext}
\usepackage[]{geometry}
\usepackage{mathrsfs}
\usepackage[all]{xy}
\usepackage{xstring}
\usepackage{slashed}

\def\N {\mathbb{N}}
\def\Z {\mathbb{Z}}

\def\R {\mathbb{R}}

\def\id{\mathrm{id}}
\def\h {\mathrm{H}}

\def\hol#1#2{\mathrm{Hol}_{#1}(#2)}
\def\quand{\quad\text{ and }\quad}
\def\quomma{\quad\text{, }\quad}

\def\ev{\mathrm{ev}}

\def\idmorph#1{#1_{dis}}
\def\hc#1{\mathrm{h}_{#1}}

\def\subset{\subseteq}

\def\nobr{~\hspace{-0.26em}}
\def\maps{\nobr:\nobr}
\def\df{\nobr := \nobr}
\def\eq{\nobr = \nobr}
\def\pr{{\mathsf{pr}}}
\def\B{\mathcal{B}\hspace{-0.03em}}
\let\Oldin\in\renewcommand{\in}{\nobr\Oldin\nobr}
\let\Oldtimes\times\renewcommand{\times}{\nobr\Oldtimes}
\let\Oldotimes\otimes\renewcommand{\otimes}{\nobr\Oldotimes}

\newlength{\widthtmp}
\def\length#1{\settowidth{\widthtmp}{#1}\the\widthtmp}

\def\lli#1{\,_{#1}\!}

\renewcommand{\varepsilon}{\epsilon}
\makeatletter
\def\bigset#1#2{\left\lbrace\;\begin{minipage}[c]{#1}\begin{center}#2\end{center}\end{minipage}\;\right\rbrace}

\makeatother

\def\erf#1{(\ref{#1})}

\newlength{\myl}

\newcommand{\ueins}{{\mathrm{U}}(1)}

\newcommand{\su}[1]{{\mathrm{SU}}\brackets{#1}}

\newcommand{\spin}[1]{{\mathrm{Spin}}\brackets{#1}}

\newcommand{\str}[1]{{\mathrm{String}}\brackets{#1}}

\def\grpd{\mathcal{G}\!rpd}

\def\trivcon#1{\mathcal{T}\!\!riv^{\nabla\!}(#1)}

\def\brackets#1{\IfStrEq{#1}{-}{}{(#1)}}
\def\buntech#1#2{\mathcal{B}\hspace{-0.01em}un_{\hspace{-0.1em}#1}^{#2}}
\def\bun#1#2{\buntech{#1}{}\brackets{#2}}

\def\bunconflat#1#2{\buntech#1{\nabla_{\!0}}\hspace{-0.05em}\brackets{#2}}

\def\grbtech#1{\mathcal{G}\hspace{-0.06em}r\hspace{-0.06em}b_{\hspace{-0.07em}#1}}
\def\grb#1#2{\grbtech#1\brackets{#2}}
\def\grbcon#1#2{\grbtech{#1}^{\nabla\!}\brackets{#2}}

\def\ugrb#1{\grb{\,}{#1}}
\def\ugrbcon#1{\grbcon\relax{#1}}

\newcommand{\alxydim}[2]{\begin{aligned}\xymatrix#1{#2}\end{aligned}}

\renewcommand{\to}{\nobr\!\xymatrix@R=0cm@C=1.4em{\ar[r] &}\nobr}
\renewcommand{\mapsto}{\!\xymatrix@R=0cm@C=1.4em{\ar@{|->}[r] &}\!}
\renewcommand{\Rightarrow}{\!\xymatrix@R=0cm@C=1.4em{\ar@{=>}[r] &}\!}
\renewcommand{\Leftarrow}{\!\xymatrix@R=0cm@C=1.4em{\ar@{<=}[r] &}\!}
\newcommand{\incl}{\!\xymatrix@R=0cm@C=1.4em{\ar@{^(->}[r] &}\!}

\renewcommand\Leftrightarrow{\!\xymatrix@R=0cm@C=1.4em{\ar@{<=>}[r] &}\!}

\usepackage{amsthm}
\usepackage{graphicx}
\usepackage[margin=2cm,font=normalsize,labelfont=bf]{caption}
\usepackage{verbatim}
\usepackage{enumerate}
\usepackage{fancyhdr}
\usepackage[]{geometry}
\usepackage[normalem]{ulem}
\usepackage{xstring}
\usepackage{needspace}

\makeatletter

\newcounter{denseversion}
\newcounter{authorcounter}
\newcounter{adresscounter}

\def\title#1{\gdef\@title{#1}}
\def\@title{}
\def\subtitle#1{\gdef\@subtitle{#1}}
\def\@subtitle{}

\def\authortagsused{0}
\def\adresstag#1{\if!#1!\else$^{\;#1\;}$\fi}

\renewcommand{\author}[2][]{
  \stepcounter{authorcounter}
  \if!#1!\else\gdef\authortagsused{1}\fi
  \ifnum\value{authorcounter}=1
    \def\@authorstringa{#2\adresstag{#1}}
    \def\@authorstringb{#2}
    \def\@authorstringc{#2\adresstag{#1}}
  \else
    \g@addto@macro\@authorstringa{\ and #2\adresstag{#1}}
    \g@addto@macro\@authorstringb{\ and #2}
    \g@addto@macro\@authorstringc{\\#2\adresstag{#1}}
  \fi}
\def\@author{\ifnum\value{denseversion}=0\@authorstringa\else\@authorstringb\fi}

\def\@adressstringa{}
\def\@adressstringb{}
\newcommand{\adress}[2][]{
  \stepcounter{adresscounter}
  \ifnum\value{adresscounter}=1
    \g@addto@macro\@adressstringa{\ifnum\authortagsused=0\def\br{\\}\else\def\br{, }\fi\adresstag{#1}#2}
    \g@addto@macro\@adressstringb{\def\br{\\}\adresstag{#1}\parbox[t]{14cm}{#2}}
  \else
    \g@addto@macro\@adressstringa{\\[\bigskipamount]\adresstag{#1}#2}
    \g@addto@macro\@adressstringb{\\[\medskipamount]\adresstag{#1}\parbox[t]{14cm}{#2}}
  \fi}
\def\@adress{\ifnum\value{denseversion}=0\@adressstringa\else\@adressstringb\fi}

\def\preprint#1{\gdef\@preprint{#1}}
\def\@preprint{}
\def\keywords#1{\gdef\@keywords{#1}}
\def\@keywords{}
\def\msc#1{\gdef\@msc{#1}}
\def\@msc{}
\def\email#1{
   \gdef\@email{#1}
   \g@addto@macro\@authorstringc{ {\it (#1)}}}
\def\@email{}
\def\dedication#1{\gdef\@dedication{#1}}
\def\@dedication{}

\def\mybaselinestretch#1{\gdef\@mybaselinestretch{#1}}
\def\@mybaselinestretch{}

\def\refname{References}

\mybaselinestretch{1.2}
\renewcommand{\baselinestretch}{\@mybaselinestretch}

\def\denseversion{
  \setcounter{denseversion}{1}
  \newgeometry{left=2cm,right=2cm,top=2cm}
  \mybaselinestretch{1.1}
  \renewcommand{\baselinestretch}{\@mybaselinestretch}
  \normalfont
  \fancyfoot[C]{\itshape{\hspace{2.5cm}--$\,\,$\thepage$\,\,$--}}}

\newlength{\myparskip}
\newlength{\myproofparskip}

\renewcommand{\emph}[1]{\def\reserved@a{it}\ifx\f@shape\reserved@a\uline{#1}\else\textit{#1}\fi}

\newcommand{\mytableofcontents}{
   \ifnum\value{denseversion}=0
     \tableofcontents
   \else
     \renewcommand{\baselinestretch}{0.8}
     \normalfont
     \tableofcontents
     \renewcommand{\baselinestretch}{\@mybaselinestretch}
     \normalfont
   \fi}

\def\href#1#2{#2}

\def\kohyp{
  \usepackage{hyperref}
  \hypersetup{
    linktocpage = true,
    pdftitle = {\@title},
    pdfauthor = {\@author},
    pdfkeywords = {\@keywords},    
    bookmarksopen = true,
    bookmarksopenlevel = 1
  }}  
\def\showkeywords{\begin{flushleft}\footnotesize\textbf{Keywords}: \@keywords\end{flushleft}}
\def\showmsc{\begin{flushleft}\footnotesize\textbf{MSC 2010}: \@msc\end{flushleft}}

\newcounter{mythm}[subsection]
\newcounter{mainthm}

\def\setsecnumdepth#1{
  \setcounter{secnumdepth}{#1}
  \setcounter{mythm}{0}
  \ifnum \c@secnumdepth >0
    \ifnum \c@secnumdepth >1
      \def\themythm{\thesubsection.\arabic{mythm}}
      \numberwithin{equation}{subsection}
      \renewcommand\theequation{\thesubsection.\arabic{equation}}
    \else
      \def\themythm{\thesection.\arabic{mythm}}
      \numberwithin{equation}{section}
      \renewcommand\theequation{\thesection.\arabic{equation}}
    \fi
  \else
    \def\themythm{\arabic{mythm}}
  \fi}

\newenvironment{mythmenv}{\strut\ \setlength{\parskip}{\myproofparskip}}{\setlength{\parskip}{\myparskip}}

\newlength{\mythmskip}
\newlength{\mythmtopskip}
\newtheoremstyle{mythmstylea}{\mythmtopskip}{\mythmskip}{\it}{}{\bf}{.}{0em}{}
\newtheoremstyle{mythmstyleb}{\mythmtopskip}{\mythmskip}{}{}{\bf}{.}{0em}{}

\theoremstyle{mythmstylea}
\newtheorem{mytheorem}[mythm]{Theorem}
\newtheorem{mydefinition}[mythm]{Definition}
\newtheorem{mycorollary}[mythm]{Corollary}
\newtheorem{myproposition}[mythm]{Proposition}
\newtheorem{mylemma}[mythm]{Lemma}

\newtheorem{mymaintheorem}[mainthm]{Theorem}
\newtheorem{mymaincorollary}[mainthm]{Corollary}
\newtheorem{mymainproposition}[mainthm]{Proposition}
\newtheorem{mymaindefinition}[mainthm]{Definition}

\theoremstyle{mythmstyleb}
\newtheorem{myremark}[mythm]{Remark}
\newtheorem{myexample}[mythm]{Example}
\newtheorem{myexercise}[mythm]{Exercise}

\newenvironment{theorem}[1][]{\begin{mytheorem}[#1]\begin{mythmenv}}{\end{mythmenv}\end{mytheorem}}
\newenvironment{definition}[1][]{\begin{mydefinition}[#1]\begin{mythmenv}}{\end{mythmenv}\end{mydefinition}}

\newenvironment{proposition}[1][]{\begin{myproposition}[#1]\begin{mythmenv}}{\end{mythmenv}\end{myproposition}}
\newenvironment{lemma}[1][]{\begin{mylemma}[#1]\begin{mythmenv}}{\end{mythmenv}\end{mylemma}}
\newenvironment{remark}[1][]{\begin{myremark}[#1]\begin{mythmenv}}{\end{mythmenv}\end{myremark}}
\newenvironment{example}[1][]{\begin{myexample}[#1]\begin{mythmenv}}{\end{mythmenv}\end{myexample}}

\renewenvironment{proof}[1][Proof]{\noindent #1. \begin{mythmenv}}{\hfill$\square$\end{mythmenv}\medskip}

\def\mytitle{}
\def\zmptitle{
  \begin{tabular}{cc}
    \begin{minipage}[c]{0.4\textwidth}
      \begin{flushleft}
        \includegraphics[width=110pt]{../../tex/zmp}
      \end{flushleft}  
    \end{minipage}&
    \begin{minipage}[c]{0.55\textwidth}
      \begin{flushright}
      {\small\sf\@preprint}
      \end{flushright}
    \end{minipage}
  \end{tabular}
  \vskip 2cm}

\def\maketitle{
  \setlength{\parskip}{\myparskip}  
  \newpage
  \noindent
  \mytitle
  \begin{center}
    \LARGE\@title\\
    \if!\@subtitle!\else \smallskip\LARGE\@subtitle\\\fi
    \bigskip
    \if!\@author!\else\bigskip\large\@author\\\fi
    \ifnum\value{denseversion}=0
      \if!\@adress!\else     \bigskip\normalsize\@adress\\\fi
      \if!\@email!\else\ifnum\value{authorcounter}=1\bigskip\normalsize\textit{\@email}\\\else\fi\fi
    \else
    \fi
    \if!\@dedication!\else \bigskip\normalsize{\@dedication}\\\fi
  \end{center}
  \ifnum\value{denseversion}=0\vskip 1.5cm\else\vskip0.5cm\fi
  \thispagestyle{empty}}

\def\kobiburl#1{
   \IfBeginWith
     {#1}
     {http://arxiv.org/abs/}
     {\kobibarxiv{#1}}
     {\kobiblink{#1}}}
\def\kobibarxiv#1{\href{#1}{\texttt{[arxiv:\StrGobbleLeft{#1}{21}]}}}
\def\kobiblink#1{Available as: \href{#1}{\texttt{#1}}}

\newcommand{\etalchar}[1]{$^{#1}$}
\def\kobib#1{
  \begin{raggedright}
  \ifnum\value{denseversion}=0\else\small\fi

  \end{raggedright}
  \ifnum\value{denseversion}=0\else
      \noindent
      \if!\@authorstringc!\else
        \ifnum\authortagsused=0\ifnum\value{authorcounter}>1\normalsize\@authorstringc\\[\medskipamount]\else\fi\else\normalsize\@authorstringc\\[\medskipamount]\fi       \fi
      \if!\@adress!\else\normalsize\@adress\\\fi
      \ifnum\authortagsused=0\ifnum\value{authorcounter}=1\if!\@email!\else\linebreak\normalsize\textit{\@email}\\\fi\else\fi\else\fi
  \fi}

\newenvironment{commentfigure}{}
\newenvironment{sidewayscommentfigure}{\begin{minipage}}{\end{minipage}}

\def\showcomments{ -- Comments suppressed}

\newif\if@fewtab\@fewtabtrue{
  \count255=\time\divide\count255 by 60
  \xdef\hourmin{\number\count255}
  \multiply\count255 by-60\advance\count255 by\time
  \xdef\hourmin{\hourmin:\ifnum\count255<10 0\fi\the\count255}}
\def\ps@draft{
  \let\@mkboth\@gobbletwo
  \def\@oddfoot{
    \hbox to 7 cm{\tiny \versionno\hfil}
    \hskip -7cm\hfil\rm\thepage\hfil{\tiny\draftdate}}
  \def\@oddhead{}
  \def\@evenhead{}
  \let\@evenfoot\@oddfoot}
\def\draftdate{\number\month/\number\day/\number\year\ \ \ \hourmin }
\newcommand\version[1]{
  \typeout{}\typeout{#1}\typeout{}
  \vskip-1.7cm \centerline{\fbox{{\normalsize\tt DRAFT -- #1 -- 
  \draftdate\showcomments}}} \vskip0.92cm}
\def\draft#1{
  \def\versionno{#1}
  \pagestyle{draft}\thispagestyle{draft}
  \gdef\@ntitle{\version\versionno \@title}
  \global\def\draftcontrol{1}}
\global\def\draftcontrol{0}

\makeatother

\usepackage[latin1]{inputenc}
\usepackage[english]{babel}

\def\quot#1{``#1''}

\def\px#1#2{P_{\!#2}#1}
\def\p{P}

\def\ev{\mathrm{ev}}

\def\hc#1{\mathrm{h}_{#1}}
\def\pcomp{\nobr\star\nobr}

\def\prev#1{\overline{#1}}

\def\un{\mathscr{R}}

\def\ufusbun#1{\mathcal{F}\!us\buntech{}{}(#1)}

\def\bundle{\mathcal{B}\!un}

\newtheorem{myconstruction}[mythm]{Construction}

\def\mtitle#1{\title{#1}}

\def\diagR{1.5cm}

\mtitle{A Construction of String 2-Group Models\\using a Transgression-Regression Technique}

\author{Konrad Waldorf}

\adress{Fakultät für Mathematik\br Universität Regensburg\br Universitätsstraße 31\br D--93053 Regensburg}

\dedication{Dedicated to Steven Rosenberg on the occasion of his 60th birthday}

\email{konrad.waldorf@mathematik.uni-regensburg.de}

\msc{Primary 22E67; secondary 53C08, 81T30, 58H05}

\kohyp

\usepackage{amstext}
\usepackage[normalem]{ulem}
\usepackage{amsbsy}

\begin{document}


\maketitle

\begin{abstract}
In this note we present a new construction of the string group that  ends optionally in two different contexts: strict diffeological 2-groups or finite-dimensional Lie 2-groups. It is canonical in the sense that no choices are involved; all the data is  written down and can be looked up (at least somewhere). The basis of our construction is the basic gerbe of Gaw\c edzki-Reis and Meinrenken. The main new insight is that under a transgression-regression procedure, the basic gerbe picks up a multiplicative structure
coming from the  Mickelsson product over the loop group. The conclusion of the construction is a relation between multiplicative gerbes and 2-group extensions  for which we use recent work of Schommer-Pries. 
\showmsc
\end{abstract}

\def\mult{\mathcal{M}\!ult}
\def\upi{\underline{\pi}}
\def\ext{\mathcal{E}\!x\!t}
\def\grpd#1{\Gamma_{\!#1}}
\def\spa{\hspace{-0.05em}}
\def\lgrpd{\mathcal{L}\spa i\spa e\spa \mathcal{G}\spa r\spa p\spa d}
\def\gbas{\mathcal{G}_{bas}}
\def\mp{\ast}

\setsecnumdepth{1}

\section{Introduction}

\label{sec:intro}

The string group $\str n$ is a topological group defined up to homotopy equivalence as the 3-connected cover of $\spin n$, for $n=3$ or $n>4$. Concrete models for $\str n$ have been provided by Stolz \cite{stolz4} and Stolz-Teichner \cite{stolz1}. In order to understand, e.g. the differential geometry of $\str n$, the so-called \quot{string geometry}, it is necessary to have models in better categories than topological groups. Its 3-connectedness implies that $\str n$ is a $K(\Z,2)$-fibration over $\spin n$,  so that it cannot be  a (finite-dimensional) Lie group.  Instead,  it allows  models in the following  contexts (in the order of appearance):
\begin{enumerate}[(i)]

\item
\label{eq:contexts:strictfrech}
Strict Fréchet Lie 2-groups \cite{baez9}.

\item
 Banach Lie 2-groups \cite{Henriques2008}.

\item
\label{eq:context:findim}
Finite-dimensional Lie 2-groups \cite{pries2}. 

\item
\label{eq:contexts:strictdiff}
Strict diffeological 2-groups \cite{schreiber2011}.

\item 
\label{eq:contexts:inffrech}
Fréchet Lie groups \cite{nikolausb}.

\end{enumerate} 
We recall that a strict  Lie 2-group is a   Lie groupoid equipped with a certain kind of monoidal structure. In the non-strict case the monoidal structure is generalized to a \quot{stacky} product. A \emph{2-group model for $\str n$} is a Lie 2-group $\Gamma$, possibly strict, Banach,  Fréchet or diffeological, equipped with a Lie 2-group homomorphism
\begin{equation}
\label{eq:stringcov}
\Gamma \to \spin n
\end{equation}
 such that the geometric realization of \erf{eq:stringcov} is a 3-connected cover.

The purpose of this note is to construct a new 2-group model for $\str n$, which can -- in the very last step -- either be chosen to live in the context \erf{eq:context:findim} of finite-dimensional Lie 2-groups, or in the context \erf{eq:contexts:strictdiff} of strict, diffeological 2-groups.
The strategy we pursue  is to reduce the problem of constructing 2-group models for $\str n$  to the construction of certain  \emph{gerbes} over $\spin n$. For the context \erf{eq:context:findim}, this reduction is possible due to  an equivalence of bicategories
\begin{equation}
\label{eq:eq1}
\bigset{9.9em}{Multiplicative, smooth bundle gerbes over $G$} 
\to
\bigset{10.1em}{Central Lie 2-group extensions of $G$ by $\B S^1$}\text{,}
\end{equation} 
which exists for any compact Lie group $G$ and reflects the fact that both bicategories are classified by $\h^4(BG,\Z)$ \cite{brylinski3,pries2}. The equivalence \erf{eq:eq1} is designed such that any multiplicative bundle gerbe over $\spin n$ whose class is a generator of $\h^4(B\spin n,\Z)\cong \Z$ automatically goes to a Lie 2-group model for $\str n $ in the context  \erf{eq:context:findim} \cite{pries2}. 
A version of the equivalence \erf{eq:eq1} exists in the  strict, diffeological context  \erf{eq:contexts:strictdiff}. 

Sections \ref{sec:multgrb} and \ref{sec:multgrb2grps} review the notions of bundle gerbes and multiplicative structures, and discuss the equivalence \erf{eq:eq1}. In Section \ref{sec:diff} we upgrade to the diffeological version. 
The following two sections are concerned with the construction of the input data, certain multiplicative bundle gerbes.  In short, the construction goes as follows: Gaw\c edzki-Reis \cite{gawedzki1,gawedzki2} and Meinrenken \cite{meinrenken1} have described a canonical construction of a bundle gerbe $\gbas$ over a compact, simple, simply-connected Lie group $G$, whose Dixmier-Douady class generates $ \h^3(G,\Z) \cong \Z$. A direct construction of a \emph{multiplicative structure} on $\gbas$ is  not known -- this is the main problem we solve in this note.

We use a transgression-regression technique developed in a series of papers \cite{waldorf9,waldorf10,waldorf11}.
The transgression of $\gbas$ is a  principal $S^1$-bundle $L\gbas$ over the loop group $LG$. Our main insight  is to combine two  additional structures one naturally finds on $L\gbas$: the \emph{Mickelsson product} \cite{mickelsson1}
and the \emph{fusion product} \cite{waldorf9}.  The fusion product allows one to \emph{regress} $L\gbas$ to a new, diffeological bundle gerbe $\un(L\gbas)$ over $G$. The Mickelsson product regresses alongside to a \emph{strictly multiplicative structure} on  $\un(L\gbas)$. Regression is inverse to transgression in the sense of a natural isomorphism
\begin{equation}
\label{eq:caniso}
\gbas \cong \un(L\gbas)
\end{equation}
of bundle gerbes over $G$. Since $\h^4(BG,\Z) \cong  \h^3(G,\Z)$ for the class of Lie groups we are looking at here, this implies that the class of $\un(L\gbas)$ generates $\h^4(BG,\Z)$.

We conclude our construction in Section \ref{sec:string} by either feeding the strictly multiplicative, diffeological bundle gerbe
$\un(L\gbas)$  into the strict, diffeological version of the equivalence  \erf{eq:eq1}, or we conclude by using the isomorphism \erf{eq:caniso} to induce a finite-dimensional,  non-strict multiplicative structure on $\gbas$ and feeding that into the equivalence \erf{eq:eq1}. For $G=\spin n$, this yields the two new 2-group models for $\str n$ in the contexts \erf{eq:contexts:strictdiff} and \erf{eq:context:findim}, respectively. 

The construction in the context \erf{eq:context:findim} is probably the most interesting result of this note. It can be seen as a small addendum to the work of Schommer-Pries \cite{pries2}. Indeed, the model of \cite{pries2}  is only  defined up to a \quot{contractible choice of isomorphisms}, while our model is canonical \quot{on the nose}.

\medskip

\noindent
{\bf Acknowledgements. }
I thank Thomas Nikolaus for discussions, in particular for  providing the argument given in Footnote \ref{foot:4}. I also thank  the Hausdorff Research Institute for Mathematics in Bonn for  kind hospitality and financial support.

\setsecnumdepth{1}
\section{Multiplicative Bundle Gerbes}

\label{sec:multgrb}

In this section we review the notion of a multiplicative bundle gerbe, which is central for this note.
Let $M$ be a smooth manifold. 

\begin{definition}
\label{def:product}
Let $\pi:Y \to M$ be a surjective submersion, and let $P$ be a principal $S^1$-bundle over the two-fold fibre product $Y^{[2]} := Y \times_M Y$. A \emph{gerbe product} on $P$ is an isomorphism
\begin{equation*}
\mu: \pr_{12}^{*}P \otimes \pr_{23}^{*}P \to \pr_{13}^{*}P
\end{equation*}
of bundles over $Y^{[3]}$ that is associative over $Y^{[4]}$.
\end{definition}

In this definition, we have denoted by $\pr_{ij}:Y^{[3]} \to Y^{[2]}$ the projection to the indexed factors, and we have denoted by $\otimes$ the tensor product of  $S^1$-bundles. Thus, a gerbe product is  for every point $(y_1,y_2,y_3) \in Y^{[3]}$ a smooth, equivariant map
\begin{equation*}
\mu: P_{(y_1,y_2)} \otimes P_{(y_2,y_3)} \to P_{(y_1,y_3)}
\end{equation*}
between  fibres of $P$. The associativity condition is that
\begin{equation*}
\mu(\mu(q_{12} \otimes q_{23}) \otimes q_{34}) = \mu(q_{12} \otimes \mu(q_{23} \otimes q_{34}))
\end{equation*}
for all $q_{ij} \in P_{(y_i,y_j)}$ and all $(y_1,y_2,y_3,y_4) \in Y^{[4]}$.

\begin{definition}[{{\cite{murray}}}]
\label{def:bundlegerbe}
A  \emph{bundle gerbe} over $M$ is a surjective submersion  $\pi\maps Y \to M$, a principal $S^1$-bundle $P$ over $Y^{[2]}$ and a gerbe product $\mu$ on $P$.
\end{definition}

Bundle gerbes over $M$ form a  bicategory $\ugrb M$ \cite{stevenson1,waldorf1}. In fact, they form a \emph{double category with companions} in the sense of \cite{Grandis,shulman1}. This means that  there are two types of 1-morphisms,  \quot{general} ones and \quot{simple} ones, together with a certain map that assigns to each simple 1-morphism a general one, its \quot{companion}.  In the case of bundle gerbes, we  call the general 1-morphisms \emph{1-isomorphisms} (they are all invertible) and the simple ones \emph{refinements}\footnote{Sometimes the simpler ones are called \quot{morphisms}, and the general ones \quot{stable isomorphisms}}. For the  definition of a 1-isomorphism we refer to \cite[Definition 5.1.2]{nikolaus}. A refinement $f: \mathcal{G} \to \mathcal{G}'$ between two bundle gerbes is  a smooth map $f_1\maps Y \to Y'$ that commutes with the two submersions to $M$, together with a bundle isomorphism $f_2: P \to P'$ over the induced map $Y^{[2]} \to Y^{\prime[2]}$, such that $f_{2}$ is a homomorphism for the gerbe products $\mu$ and $\mu'$. The assignment of a 1-morphism to a refinement can be found in \cite[Lemma 5.2.3]{nikolaus}.

The bicategory $\ugrb M$ is equipped with many additional features. For instance, it is monoidal, and the assignment $M \mapsto \ugrb M$ is a  sheaf of monoidal bicategories over the site of smooth manifolds (with surjective submersions) \cite{stevenson1,waldorf1,nikolaus}. This means in particular that one can consistently pull back and tensor bundle gerbes, refinements, 1-isomorphisms, and 2-morphisms. 
Denoting by  $\hc 0$ the operation of taking the set of isomorphism classes we have:

\begin{theorem}[{{\cite{murray2}}}]
\label{th:grbclass}
$\hc 0 \ugrb  M \cong \h^3(M,\Z)$.
\end{theorem}

In the following we consider bundle gerbes over a \emph{Lie group} $G$. For preparation, let us suppose that $\pi:Y \to G$ is a surjective submersion, such that $Y$ is another Lie group and $\pi$ is a group homomorphism. Then, the fibre products $Y^{[k]}$ are again Lie groups, and the projections $\pr_{ij} \maps Y^{[3]} \to Y^{[2]}$ are Lie group homomorphisms. Suppose further that we have a \emph{central extension}
\begin{equation*}
1 \to S^1 \to P \to Y^{[2]} \to 1
\end{equation*}
of Lie groups, i.e. a central extension of groups such that $P$ is a principal $S^1$-bundle over $Y^{[2]}$. In this situation, a  gerbe product $\mu$ on $P$ is called \emph{multiplicative} if it is a group homomorphism, i.e. if
\begin{equation*}
\mu(p_{12}^{}p_{12}' \otimes p_{23}^{}p_{23}') = \mu(p_{12} \otimes p_{23}) \cdot \mu(p_{12}' \otimes p_{23}')
\end{equation*}   
for all $p_{ij} \in P_{(y_i,y_j)}$, $p_{ij}'\in P_{(y_i',y_j')}$  and all $(y_1,y_2,y_3),(y_1',y_2',y_3')\in Y^{[3]}$.

\begin{definition}
Let $\mathcal{G}=(Y,\pi,P,\mu)$ be a bundle gerbe over $G$. A \emph{strictly multiplicative structure} on $\mathcal{G}$ is a Lie group structure on $Y$ such that $\pi$ is a group homomorphism, together with a Lie group structure on $P$, such that $P$ is a central extension of $Y^{[2]}$ by $S^1$ and $\mu$ is multiplicative.
\end{definition}

A bundle gerbe $\mathcal{G}$ together with a strictly multiplicative structure is called a \emph{strictly multiplicative bundle gerbe}.
The problem is that strictly multiplicative structures on bundle gerbes rarely exist. The following definition is a suitable generalization.

\begin{definition}[{{\cite{brylinski3,carey4,waldorf5}}}]
\label{def:mult}
A \emph{multiplicative structure} on a bundle gerbe  $\mathcal{G}$ over $G$ is a 1-isomorphism
\begin{equation*}
\mathcal{M}: \pr_1^{*}\mathcal{G} \otimes \pr_2^{*}\mathcal{G}
\to m^{*}\mathcal{G}
\end{equation*}
of bundle gerbes over $G \times G$, and a 2-isomorphism
\begin{equation*}
\alxydim{@R=\diagR@C=2.1cm}{\mathcal{G}_1 \otimes \mathcal{G}_2 \otimes \mathcal{G}_3 \ar[r]^-{\mathcal{M}_{1,2} \otimes \id} \ar[d]_{\id \otimes \mathcal{M}_{2,3}} & \mathcal{G}_{12} \otimes \mathcal{G}_3  \ar[d]^{\mathcal{M}_{12,3} } \ar@{=>}[dl]|*+{\alpha} \\ \mathcal{G}_1 \otimes \mathcal{G}_{23}  \ar[r]_-{\mathcal{M}_{1,23} } & \mathcal{G}_{123} }
\end{equation*}
between 1-isomorphisms over $G \times G \times G$ that satisfies
the obvious   pentagon axiom.
\end{definition}

In this definition, $m: G \times G \to G$ denotes the multiplication of $G$, and the index convention is such that e.g. the index $(..)_{ij,k}$ stands for the pullback along the map $(g_i,g_j,g_k) \mapsto (g_ig_j,g_k)$. For instance, $\mathcal{G}_i = \pr_i^{*}\mathcal{G}$ and $\mathcal{G}_{12}=m^{*}\mathcal{G}$. A \emph{multiplicative bundle gerbe} over $G$ is a bundle gerbe together with a multiplicative structure. Multiplicative bundle gerbes over $G$ form a  bicategory that we denote by $\mult\ugrb G$. We have for compact Lie groups $G$:

\begin{theorem}[{{\cite[Propositions 1.5 and 1.7]{brylinski3}}}]
\label{th:multgrbclass}
$\hc 0 \mult\ugrb G \cong  \h^4(BG,\Z)$.
\end{theorem}

A \emph{strictly} multiplicative structure on a bundle gerbe $\mathcal{G}=(Y,\pi,P,\mu)$ induces a multiplicative structure in the following way. 
Over $G \times G$, we consider the bundle gerbes $\mathcal{G}_{1,2} = \pr_1^{*}\mathcal{G} \otimes \pr_2^{*}\mathcal{G}$ and $\mathcal{G}_{12} = m^{*}\mathcal{G}$. Employing the definitions of pullbacks and tensor products \cite{waldorf1}, the bundle gerbe $\mathcal{G}_{1,2}$ consists of $Y_{1,2} := Y \times Y$ with   $\pi_{1,2} :=\pi \times \pi$, and the  bundle $P^{1,2}$ over  $Y_{1,2}^{[2]}$ with fibres
\begin{equation*}
P^{1,2}_{(y_1,y_2),(y_1',y_2')} = P_{y_1,y_1'} \otimes P_{y_2,y_2'}\text{.}
\end{equation*}
The bundle gerbe product $\mu^{1,2}$ on $P^{1,2}$ is just the tensor product of $\mu$ with itself.
On the other side, the bundle gerbe $\mathcal{G}_{12}$ is $Y_{12} := G \times Y$ with  $\pi_{12}(g,y) := (g,g^{-1}\pi(y))$,   while the principal bundle $P^{12}$ and the bundle gerbe product $\mu^{12}$ are just pullbacks along the projection $Y_{12} \to Y$. 
Now, the multiplication of the Lie group $Y$ defines a smooth map
\begin{equation*}
f_1: Y_{1,2} \to Y_{12}: (y_1,y_2) \mapsto (\pi(y_1),y_1y_2)
\end{equation*}
that commutes with  $\pi_{1,2}$ and $\pi_{12}$. Further, the multiplication on $P$ defines a bundle isomorphism
\begin{equation*}
f_2: P^{1,2} \to P^{12}: (p,p')\mapsto pp'
\end{equation*}
and  the multiplicativity of $\mu$ assures that $f_2$ is a homomorphism for the bundle gerbe products $\mu^{1,2}$ and $\mu^{12}$. Thus, the pair $(f_1,f_2)$ is a refinement $f: \mathcal{G}_{1,2} \to \mathcal{G}_{12}$ which in turn defines the required 1-isomorphism $\mathcal{M}$. 
Next we look at the diagram over $G \times G \times G$ of Definition \ref{def:mult}. It turns out that the associativity of the Lie groups $Y$ and $P$ imply  the strict commutativity of the refinements representing the four 1-isomorphisms in the diagram. In this case, the coherence of companions in double categories provides the required 2-isomorphism $\alpha$, and a general coherence result implies the  pentagon axiom. This concludes the construction of a multiplicative bundle gerbe $(\mathcal{G},\mathcal{M},\alpha)$ from a strictly multiplicative one.

The 2-functor $\mult\ugrb G \to \ugrb G$ that forgets the multiplicative structure corresponds \cite[Lemma 2.3.9]{waldorf5} under the  bijections of Theorems \ref{th:grbclass} and \ref{th:multgrbclass} to   the usual  \quot{transgression} map
\begin{equation}
\label{eq:transZ}
\h^4(BG,\Z) \to \h^3(G,\Z)\text{.}
\end{equation}
If $G$ is compact, simple, and simply connected, this map is a bijection, so that every bundle gerbe over $G$ has a (up to isomorphism) unique multiplicative structure. 
If $G$ is  only compact and simple,  the map \erf{eq:transZ} is still injective, but the existence of multiplicative structures is obstructed\footnote{This can e.g. be seen  by looking at  the descent theory for multiplicative  gerbes  \cite{gawedzki9}.}.

\setsecnumdepth{2}

\section{Lie 2-Group Extensions}

\label{sec:multgrb2grps}

We relate multiplicative bundle gerbes to central Lie 2-group extensions. The material presented here is well-known; the whole section can be seen as an expansion of \cite[Remark 101]{pries2}.

\subsection{Lie 2-Groups}
 
\label{sec:lie2groups} 
 
We recall that a \emph{Lie groupoid} is a groupoid $\Gamma$ whose objects $\Gamma_0$ and morphisms $\Gamma_1$ form smooth manifolds, whose source and target maps are surjective submersions, and whose composition and inversion are smooth maps. 

\begin{example} 
\label{ex:liegroupoids}
\begin{enumerate}[(i)]

\item 
\label{Xdis}
Every smooth manifold $X$ defines a \quot{discrete} Lie groupoid $\idmorph{X}$ with objects $X$ and only identity morphisms.

\item
\label{BG}
Every Lie group $G$ defines a  Lie groupoid $\B G$ with one object and automorphism  group  $G$. 

\item
\label{ex:gerbe}
Let $\mathcal{G} = (Y,\pi,P,\mu)$ be a bundle gerbe over $M$. Then, we have a Lie groupoid
 $\grpd{\mathcal{G}}$ with objects  $Y$ and morphisms  $P$. Source and target maps are defined by $s := \pr_1 \circ \chi$ and $t := \pr_2 \circ \chi$, where $\chi:P \to Y^{[2]}$ denotes the bundle projection, and the composition is the gerbe product $\mu$. Identities and inversion are also induced by $\mu$ \cite[Corollary 5.2.6 (iii)]{nikolaus}.

\end{enumerate}
\end{example}

Lie groupoids -- like bundle gerbes -- form a double category with companions, denoted  $\lgrpd$. The simple 1-morphisms are smooth functors. The  general ones are   \emph{smooth anafunctors}\footnote{Sometimes smooth anafunctors are called   \quot{Hilsum-Skandalis morphisms} or \quot{bibundles}.}  $P: \Gamma \to \Omega$, which are principal $\Omega$-bundles $P$ over $\Gamma$, see \cite[Section 2]{nikolaus} for a detailed discussion and references. The 2-morphisms are $\Omega$-bundle isomorphisms over $\Gamma$, and will be called \emph{smooth transformations}. 

\begin{proposition}
\label{prop:framed2functor}
The assignment $\mathcal{G}\mapsto \grpd{\mathcal{G}}$ of a Lie groupoid to a bundle gerbe extends to a 2-functor
$\ugrb M \to \lgrpd$ that respects companions.
\end{proposition}

\begin{proof}
The 2-functor is constructed in \cite[Section 7.2]{nikolaus}. The claim that this 2-functor respects companions means additionally that it sends a refinement $f: \mathcal{G} \to \mathcal{G}'$   to  a smooth functor $\grpd{\mathcal{G}} \to \grpd{\mathcal{G}'}$; this can easily be checked using the given definitions.    
\end{proof}

\begin{definition}[{{\cite{baez5}}}]
A \emph{strict Lie 2-group} is a Lie groupoid $\Gamma$ whose objects $\Gamma_0$ and morphisms $\Gamma_1$ form Lie groups, such that source, target, and  composition are group homomorphisms. 
\end{definition}

\noindent
Continuing Example \ref{ex:liegroupoids}, it is easy to check the following statements:
\begin{enumerate}[(i)]
\item 
If $G$ is a Lie group, the Lie groupoid  $\idmorph{G}$ is a strict Lie 2-group.

\item
If $A$ is an abelian Lie group, the Lie groupoid $\B A$ is a strict Lie group.

\item
If $\mathcal{G}$ is a strictly multiplicative bundle gerbe over $G$, the Lie groupoid $\grpd{\mathcal{G}}$ is a strict Lie 2-group. 
\end{enumerate}
In order to  include non-strictly multiplicative bundle gerbes, we need the following generalization:

\begin{definition}[{{\cite{baez5,pries2}}}]
\label{def:Lie2group}
A \emph{Lie 2-group} is a Lie groupoid $\Gamma$  with  smooth anafunctors \begin{equation*}
m: \Gamma \times \Gamma \to \Gamma
\quand
e:1 \to \Gamma\text{,}
\end{equation*}
and smooth transformations $\alpha$, $l$, $r$,  where $\alpha$ expresses that $m$ is an associative product and $l,r$ express that $e$ is a left and right unit for this product, such that the smooth anafunctor
\begin{equation*}
(\pr_1,m): \Gamma \times \Gamma \to \Gamma \times \Gamma
\end{equation*}
is invertible.
\end{definition}

In this definition, $1$ denotes the trivial Lie groupoid. The details about the smooth transformations can e.g. be found in \cite[Definition 41]{pries2}. 
We have the following examples of Lie 2-groups:

1.) If $\Gamma$ is a \emph{strict Lie 2-group}, the Lie group structures on $\Gamma_0$ and $\Gamma_1$ can be bundled into smooth functors $m: \Gamma \times \Gamma \to \Gamma$ and $e: 1 \to \Gamma$  satisfying strictly the axioms of an associative multiplication and of a unit. The coherence  of companions in the double category  $\lgrpd$ provides associated smooth anafunctors and the required  smooth transformations. Thus,  \emph{strict} Lie 2-groups are particular Lie 2-groups.  

2.) A multiplicative structure $(\mathcal{M},\alpha)$ on a bundle gerbe $\mathcal{G}$ equips the Lie groupoid $\grpd{\mathcal{G}}$  with a Lie 2-group structure. Indeed, one can check explicitly that $\grpd{\pr_1^{*}\mathcal{G} \otimes \pr_2^{*}\mathcal{G}} = \grpd{\mathcal{G}} \times \grpd{\mathcal{G}}$ as Lie groupoids, and also produce an evident smooth functor $\pr : \grpd{m^{*}\mathcal{G}} \to \grpd{\mathcal{G}}$. Using that $\Gamma$ is functorial (Proposition \ref{prop:framed2functor}), we obtain a smooth anafunctor
\begin{equation*}
\alxydim{}{\grpd{\mathcal{G}} \times \grpd{\mathcal{G}} = \grpd{\pr_1^{*}\mathcal{G} \otimes \pr_2^{*}\mathcal{G}} \ar[r]^-{\grpd{\mathcal{M}}} & \grpd{m^{*}\mathcal{G}} \ar[r]^{\pr} & \grpd{\mathcal{G}}\text{.} }
\end{equation*} 
Similarly, one can check that the 2-isomorphism $\alpha$ provides the required associator for this multiplication. Let $1: pt \to G$ denote the unit element of the group $G$. Using duals of bundle gerbes one can show that the 1-isomorphism $\mathcal{M}$ induces  a distinguished 1-isomorphism $\mathcal{E}: \mathcal{I} \to 1^{*}\mathcal{G}$, where $\mathcal{I}$ is the trivial $S^1$-bundle gerbe over the point. We have $\grpd{\mathcal{I}} = \B S^1$, and obtain, again by functorality of $\Gamma$, the required smooth anafunctor
\begin{equation*}
\alxydim{}{1 \ar[r] & \B S^1 \ar[r]^{\grpd{\mathcal{E}}} & \grpd{1^{*}\mathcal{G}} \ar[r]^{\pr} & \grpd{\mathcal{G}}\text{.}}
\end{equation*}
The smooth transformations $l$ and $r$ can both be deduced from the 2-isomorphism $\alpha$. 

\subsection{Central Extensions}

\label{sec:extensions}

We briefly review some aspects of  principal 2-bundles \cite{bartels,pries2}.  Let $\Gamma$ be a   Lie 2-group.   A \emph{principal $\Gamma$-2-bundle} over a smooth manifold $M$ is a Lie groupoid $\mathcal{P}$ \quot{total space}, a smooth functor $\pi: \mathcal{P} \to \idmorph{M}$ \quot{projection}, a smooth anafunctor
$\tau: \mathcal{P} \times \Gamma \to \mathcal{P}$ \quot{right action} together with two smooth transformations satisfying several axioms. 
If $\Gamma$ is a strict Lie 2-group, a principal $\Gamma$-2-bundle is called \emph{strict} if  $\tau$ is a smooth \emph{functor}, and both smooth transformations are identities. Strict principal $\Gamma$-2-bundles have been studied  in detail in \cite{nikolaus}.

\begin{example}
\label{ex:2bundle}
We recall from Example \ref{ex:liegroupoids} \erf{ex:gerbe}  that there is a Lie groupoid $\grpd{\mathcal{G}}$ associated to any bundle gerbe $\mathcal{G}$ over $M$. Together with the smooth functor  $\pi\maps \grpd{\mathcal{G}} \to \idmorph{M}$ given by the surjective submersion of $\mathcal{G}$, and the smooth functor $\tau: \grpd{\mathcal{G}} \times \B S^1 \to \grpd{\mathcal{G}}$ induced by the   action  of $S^1$ on $P$, this yields a strict principal $\B S^1$-2-bundle over $M$, see \cite[Example 73]{pries2}, \cite[Section 7.2]{nikolaus}.
\end{example}

\begin{proposition}[{{\cite[Theorem 7.1]{nikolaus}}}]
\label{prop:stackeq}
Example \ref{ex:2bundle} establishes an equivalence between the bicategories of bundle gerbes over $M$ and strict principal $\B S^1$-2-bundles over $M$.
\end{proposition}

Schommer-Pries has introduced a very general notion of Lie 2-group extensions \cite[Definition 75]{pries2}. For the purpose of this note we may reduce it to the case that a \quot{discrete} Lie 2-group $\idmorph{G}$ is extended by the \quot{codiscrete} Lie 2-group $\B S^1$.

\begin{definition}
\label{def:extension}
Let $G$ be a Lie group. 
A \emph{Lie 2-group extension} of $G$ by $\B S^1$ is a Lie 2-group $\Gamma$ with Lie 2-group homomorphisms
\begin{equation*}
\alxydim{}{\B S^1 \ar[r]^-{i} & \Gamma \ar[r]^-{\pi} & \idmorph{G}}
\end{equation*}
such that:
\begin{enumerate}[(i)]
\item 
The composite $\pi \circ i$ is  the constant functor $1: \B S^1 \to \idmorph{G}$.

\item
$\pi:\Gamma \to \idmorph{G}$ is a principal $\B S^1$-2-bundle over $G$.
\end{enumerate}
The extension is called \emph{strict} if $\Gamma$ is a strict Lie 2-group, and  $\pi$, $i$ are  strict  2-group homomorphisms.
\end{definition}

The notion of \emph{central} Lie 2-group extensions introduced in \cite[Definition 83]{pries2} requires a certain group homomorphism $\alpha:G \to \mathrm{Aut}(S^1) \cong \Z/2\Z$ to be trivial. Central Lie 2-group extensions of $G$ by $\B S^1$ form a bicategory $\ext (G,\B S^1)$, and for $G$ compact  we have:

\begin{theorem}[{{\cite{pries2}}}]
\label{th:extclass}
 $\hc 0 \ext (G,\B S^1) \cong \h^4(BG,\Z)$.
\end{theorem}

As discussed in Section \ref{sec:lie2groups}, the Lie groupoid $\grpd{\mathcal{G}}$ associated to a strictly multiplicative bundle gerbe  $\mathcal{G}=(Y,\pi,P,\mu)$  over $G$ is a strict Lie 2-group. We have  the functor $\pi\maps \grpd{\mathcal{G}} \to \idmorph{G}$ from Example \ref{ex:2bundle}, and a functor $i: \B S^1 \to \grpd{\mathcal{G}}$ induced by the second arrow of the central
extension
\begin{equation}
\label{eq:centexrep}
1 \to S^1 \to P \to Y^{[2]} \to 1\text{.}
\end{equation}
Condition (i) is clear, and (ii) is  proved by Example \ref{ex:2bundle}. Centrality follows from the one of \erf{eq:centexrep}. Thus, every strictly multiplicative bundle gerbe defines a central, strict Lie 2-group extension.

If $\mathcal{G}$ is a multiplicative bundle gerbe over $G$, the Lie groupoid  $\grpd{\mathcal{G}}$ is a Lie 2-group. The functor   $\pi: \grpd{\mathcal{G}} \to \idmorph{G}$ is the same as before, and the smooth anafunctor $i: \B S^1 \to \grpd{\mathcal{G}}$ is defined by
\begin{equation*}
\alxydim{@=1cm}{\mathcal{B}S^1 = pt \times \mathcal{B}S^1 \ar[r]^-{e \times \id} & \grpd{\mathcal{G}} \times \mathcal{B}S^1 \ar[r]^-{\tau} & \grpd{\mathcal{G}}\text{,}}
\end{equation*}
where $\tau$ is the action functor of Example \ref{ex:2bundle}. Conditions (i) and (ii) are still satisfied, and centrality can be concluded from the strict case, since it only affects the underlying \quot{discrete} 2-groups and every Lie 2-group can be strictified upon discretization. Thus, every multiplicative bundle gerbe defines a central Lie 2-group extension. 
Summarizing, we obtain the following  (commutative) diagram of bicategories and 2-functors:
\begin{equation}
\label{eq:diagram}
\alxydim{@C=\diagR@R=\diagR}{\left \lbrace \txt{Strictly multiplicative\\bundle gerbes over $G$} \right\rbrace \ar[d] \ar[r] & \left \lbrace\txt{Central strict Lie 2-group\\extensions of $G$ by $\B S^1$}\right\rbrace \ar[d] \\ \left \lbrace\txt{Multiplicative bundle\\gerbes over $G$}\right\rbrace \ar[r] & \left \lbrace\txt{ Central Lie 2-group\\extensions of $G$ by $\B S^1$}\right\rbrace}
\end{equation}

\begin{theorem}
\label{th:multgrbext}
The horizontal 2-functors in  diagram \erf{eq:diagram} are equivalences of bicategories.
If $G$ is compact, they induce the identity on $\h^4(BG,\Z)$ under the bijections of Theorems \ref{th:multgrbclass} and \ref{th:extclass}. \end{theorem}

\begin{proof}
For the purposes of this note, it suffices to prove the second statement for $G$ compact and simple. Then, since $\h^4(BG,\Z) \to \h^3(G,\Z)$ is injective, it suffices to observe that the horizontal 2-functors induce the identity on $\h^3(G,\Z)$. The maps to $\h^3(G,\Z)$ induced by the bijections of Theorems \ref{th:multgrbclass} and \ref{th:extclass} are, respectively,  the projection to the underlying bundle gerbe, see \erf{eq:transZ}, and the projection to the underlying principal $\B S^1$-2-bundle of a 2-group extension, see \cite{pries2}. Under both horizontal 2-functors, these are related by the assignment of Example \ref{ex:2bundle}, which is an equivalence of bicategories (Proposition \ref{prop:stackeq}).
\end{proof}

\setsecnumdepth{1}

\section{The Site of Diffeological Spaces}

\label{sec:diff}

We recall that a \emph{site} is a category together with a \emph{Grothendieck (pre-)topology}: a class of morphisms called \emph{coverings}, containing all identities, closed under composition, and  stable under  pullbacks along arbitrary morphisms.
Above we have presented the definitions of bundle gerbes, groupoids, 2-groups, and 2-group extensions internal to the familiar site $C^{\infty}$ of smooth (finite-dimensional) manifolds, with the coverings given by   surjective submersions.

Schommer-Pries  proved that the site $C^{\infty}$ allows  2-group models for the string group \cite[Theorem 2]{pries2}. However, one can show that it does not allow  \emph{strict} 2-group models\footnote{\label{foot:4}By Theorem \ref{th:multgrbext} such a strict Lie 2-group extension of $G$ by $\B S^1$ would correspond to a strictly multiplicative bundle gerbe $\mathcal{G}$ over $G$ whose Dixmier-Douady class generates $\h^3(G,\Z)$. We may assume that $G=\su 2$, otherwise we consider the restriction of $\mathcal{G}$ to an  $\su 2$ subgroup (with still non-trivial Dixmier-Douady class). The strict Lie 2-group  $\grpd{\mathcal{G}}$ induces an exact sequence
\begin{equation*}
\alxydim{@C=0.6cm}{1 \ar[r] & S^1 \ar[r] & \mathrm{ker}(s) \ar[r]^-{t} & Y \ar[r]^-{\pi} & \su 2 \ar[r] & 1}
\end{equation*}
of Lie groups \cite[Section 3]{nikolausa}, where $s,t$ are the source and target  maps of $\grpd{\mathcal{G}}$.  Thus, the submersion
 $\pi:Y \to \su 2$ of $\mathcal{G}$ is a principal bundle for the structure group $H := \mathrm{ker}(s)/S^1$. Such bundles are classified by $\pi_2(H)=0$, which implies that it has a global section, in contradiction to the non-triviality of $\mathcal{G}$, see \cite[Lemma 3.2.3]{waldorf10}.}. As mentioned in Section \ref{sec:intro}, 
strictness can be achieved by passing to a bigger site, e.g. the site $F^{\infty}$ of (possibly infinite-dimensional) Fréchet manifolds   \cite{baez9}. For the transgression-regression technique we want to use in the next section we have to pass  to a yet bigger site, the site $D^{\infty}$ of diffeological spaces.

We refer to \cite[Appendix A.1]{waldorf9} for an introduction to diffeological spaces and  references. In short, a \emph{diffeological space} is a set $X$ together with a collection of generalized charts called \quot{plots}. A \emph{plot} is a triple $(n,U,c)$ consisting of a  number $n\in \N$, an open subset $U\subset \R^n$ and a map $c: U \to X$.  A map $f:X \to X'$ between diffeological spaces is \emph{smooth} if its composition $f \circ c$ with every plot $c$ of $X$ is a plot of $X'$. This defines the category $D^{\infty}$ of diffeological spaces. A Grothendieck topology on $D^{\infty}$ is provided by \emph{subductions}: smooth maps $\pi:Y \to X$ such that every plot $c:U \to X$ lifts locally to $Y$.

A manifold $M$ (either smooth or Fréchet) can be regarded as a diffeological space with the underlying set $M$, and the plots given by \emph{all} smooth maps $c:U \to M$, for all open subsets of $\R^n$ and all $n$. We obtain a sequence
\begin{equation}
\label{eq:sites}
C^{\infty} \to F^{\infty} \to D^{\infty}
\end{equation}
of functors. These preserve the Grothendieck topologies in the sense that they send surjective submersions to subductions. Furthermore, they are full and faithful: this means that upon embedding two objects into a bigger site, the set of all maps between them is not getting bigger or smaller.

If some definition  is given in terms of the ingredients of a certain site, the same definition can obviously be repeated in any other site. For example, a \emph{smooth principal $S^1$-bundle} over a smooth manifold $X$ can be defined as a surjective submersion $\pi:P \to X$ together with a smooth map $\tau: P \times S^1 \to P$ that defines a fibrewise action, such that $(\pr_1,\tau)\maps P \nobr\times\nobr S \to P \times_X P$ is a diffeomorphism. Accordingly, a \emph{diffeological principal $S^1$-bundle} over a diffeological space $X$ is a subduction $\pi:P \to X$ and a smooth map $\tau$ satisfying the same conditions; see \cite{waldorf9} for a thorough discussion. Similarly, one repeats the definition of a bundle gerbe, of a Lie groupoid, of a Lie 2-group, and of a Lie 2-group extension in the site of diffeological spaces. 

The classification Theorems \ref{th:grbclass}, \ref{th:multgrbclass} and \ref{th:extclass} remain true for (multiplicative) diffeological bundle gerbes and diffeological 2-group extensions, as long as the base spaces $M$ and $G$ are finite-dimensional smooth manifolds; see e.g. \cite[Theorem 3.1.3]{waldorf10}. Similarly, the relation between multiplicative bundle gerbes and 2-group extensions of Theorem \ref{th:multgrbext} remains true in the diffeological context. In particular, there is a 2-functor
\begin{equation}
\label{eq:cextstrictdiff}
\bigset{14.7em}{Strictly multiplicative diffeological bundle gerbes over $G$}
\to
\bigset{13.8em}{Central, strict diffeological 2-group extensions of $G$ by $\B S^1$}
\end{equation}
that induces the identity on $\h^4(BG,\Z)$ for $G$  compact and simple. This will be used in Section \ref{sec:string}.

\section{The Transgression-Regression Machine}

\label{sec:transgressionregression}

Brylinski and McLaughlin have defined a procedure to transform a bundle gerbe over a smooth manifold $M$ into a Fréchet principal $S^1$-bundle over the Fréchet manifold $LM \df C^{\infty}(S^1,M)$, setting up an important relation between geometry on a manifold and geometry on its loop space \cite{brylinski1, brylinski4}. Their procedure uses, as an auxiliary datum, a connection on the bundle gerbe. 
If $\ugrbcon M$ denotes the bicategory of bundle gerbes with connection over $M$, and $\hc 1$ denotes the operation of producing a category (by identifying 2-isomorphic morphisms), then Brylinski's and McLaughlin's construction furnishes a functor
\begin{equation*}
L : \hc 1 \ugrbcon M \to \bun {S^1}{LM}\text{.}
\end{equation*}

We shall describe some details of the construction following \cite{waldorf5,waldorf10}. 
If $\mathcal{G}$ is a bundle gerbe with connection over $M$, the fibre of $L\mathcal{G}$ over a loop $\tau \in LM$ is
\begin{equation}
\label{eq:lgfibre}
L\mathcal{G}|_{\tau} := \hc 0 \trivcon{\tau^{*}\mathcal{G}}\text{,}
\end{equation} 
i.e. it consists of isomorphism classes of (connection-preserving) trivializations of $\tau^{*}\mathcal{G}$. In general, trivializations of a bundle gerbe $\mathcal{K}$ over a smooth manifold $X$ form a category that is a torsor for the monoidal category $\bunconflat {{S^1}} X$ of flat principal $S^1$-bundles over $X$, under a certain action functor
\begin{equation}
\label{eq:torsor}
\trivcon {\mathcal{K}} \times \bunconflat {{S^1}} X \to  \trivcon {\mathcal{K}}: (\mathcal{T},P) \mapsto \mathcal{T} \otimes P\text{.}
\end{equation}
The fibres \erf{eq:lgfibre} are thus torsors over the group $\hc 0 \bunconflat {{S^1}}{S^1} \cong S^1$. There exists a unique Fréchet manifold structure on $L\mathcal{G}$ turning it into a Fréchet principal $S^1$-bundle \cite[Proposition 3.1.2]{waldorf5}.

It is   easier to pass to the site of diffeological spaces. The plots of  $LM$ are  maps $c \maps U \to LM$ whose \emph{adjoint map} $U \times S^1 \to M: (u,z) \mapsto c(u)(z)$ is smooth (in the ordinary sense) \cite[Lemma A.1.7]{waldorf9}. The plots of the total space $L\mathcal{G}$ are  maps  $c:U \to L\mathcal{G}$ for which every point $w\in U$ has an open neighborhood $w \in W \subset U$ such that
\begin{enumerate}[(i)]

\item 
the map $d: W \times S^1 \to M$ defined by
\begin{equation*}
\alxydim{@C=1.5cm}{W \times S^1 \ar[r]^-{c|_W \times \id} & L\mathcal{G} \times S^1 \ar[r]^-{\mathrm{pr} \times \id} & LM \times S^1 \ar[r]^-{\ev} & M}
\end{equation*}
is a smooth map, and

\item
there exists a trivialization $\mathcal{T}$ of $d^{*}\mathcal{G}$ with $c(x) \cong \iota_x^{*}\mathcal{T}$ for all $x\in W$, where  $\iota_x: S^1 \to W \times S^1$ is  $\iota_x(z):=(x,z)$. 

\end{enumerate}

Diffeological principal $S^1$-bundles over $LM$  in the image of the transgression functor $L$ are  equipped with more structure.  Relevant for this note is a \emph{fusion product} \cite{waldorf10}. We denote by $PM$ the set of smooth paths $\gamma: [0,1] \to M$ with sitting instants, i.e. $\gamma$ is constant near the endpoints. This ensures that two  paths $\gamma_1$, $\gamma_2$ with a common end  can be composed to another smooth path $\gamma_2 \pcomp \gamma_1$. The set $PM$ is \emph{not} a Fréchet manifold, but a nice diffeological space whose plots are again those  maps $c:U \to PM$ whose adjoint map $(u,t) \mapsto c(u)(t)$ is smooth. The evaluation map\begin{equation*}
\ev: \p M \to M \times M: \gamma \mapsto (\gamma(0),\gamma(1))
\end{equation*}
is obviously smooth, and a subduction if $M$ is connected.  
We denote by $PM^{[k]}$ the fibre product of $PM$ over $M\times M$; it consists of $k$-tuples of  paths with a common initial point and a common end point. If we denote by $\prev\gamma$ the inverse of a path $\gamma$, we obtain a  smooth map \cite[Section 2.2]{waldorf9}
\begin{equation*}
\ell: \p M^{[2]} \to LM: (\gamma_1,\gamma_2) \mapsto \prev{\gamma_2} \pcomp \gamma_1\text{.}
\end{equation*}

\begin{definition}[{{\cite[Definition 2.1.3]{waldorf10}}}]
\label{def:fusionproduct}
Let $P$ be a diffeological principal $S^1$-bundle over $LM$. A \emph{fusion product} on $P$ is a gerbe product $\lambda$ on $\ell^{*}P$ in the sense of Definition \ref{def:product}. 
\end{definition}

Explicitly, a fusion product $\lambda$ provides, for each triple $(\gamma_1,\gamma_2,\gamma_3) \in PM^{[3]}$ a smooth map
\begin{equation*}
\lambda: P_{\prev{\gamma_2} \pcomp \gamma_1} \otimes P_{\prev{\gamma_3} \pcomp \gamma_2} \to P_{\prev{\gamma_3} \pcomp \gamma_1}\text{,}
\end{equation*}
and these maps are associative over quadruples of paths. A pair $(P,\lambda)$ is called a \emph{fusion bundle}. We  denote by $\ufusbun {LM}$ the category of fusion bundles over $LM$. 
The important point established in \cite[Section 4.2]{waldorf10} is that the  functor $L$ lifts to a functor
\begin{equation*}
\ugrbcon M \to \ufusbun {LM}\text{,}
\end{equation*}
i.e. a principal $S^1$-bundle in the image of transgression is equipped with  a \emph{canonical fusion product} $\lambda_{\mathcal{G}}$. Let us briefly recall how $\lambda_{\mathcal{G}}$ is characterized. We denote by $\iota_1,\iota_2: [0,1] \to S^1$ the inclusion of the interval into the left and the right half of the circle. Let $(\gamma_1,\gamma_2,\gamma_3)$ be a triple of paths with a common initial point $x$ and a common end point $y$, and let $\mathcal{T}_{ij}$ be trivializations of the pullback of $\mathcal{G}$ to the loops $\ell(\gamma_i,\gamma_j)$, for  $(ij)=(12),(23),(13)$. Then,
the relation\begin{equation*}
\lambda_{\mathcal{G}}(\mathcal{T}_{12} \otimes \mathcal{T}_{23}) = \mathcal{T}_{13}
\end{equation*}
holds if and only if there exist 2-isomorphisms
\begin{equation*}
\phi_1: \iota_1^{*}\mathcal{T}_{12} \Rightarrow \iota_1^{*}\mathcal{T}_{13}
\quomma
\phi_2: \iota_2^{*}\mathcal{T}_{12} \Rightarrow \iota_1^{*}\mathcal{T}_{23}
\quand
\phi_3: \iota_2^{*}\mathcal{T}_{23} \Rightarrow \iota_2^{*}\mathcal{T}_{13}
\end{equation*}
between trivializations of the pullbacks of $\mathcal{G}$ to the paths $\gamma_1$, $\gamma_2$, and $\gamma_3$, respectively, such that their restrictions to the two common points $x$ and $y$ satisfy the cocycle condition
$\phi_1 = \phi_3 \circ \phi_2$.

A fusion product permits one to define a functor inverse to transgression \cite[Section 5.1]{waldorf10}. 
Suppose $(P,\lambda)$ is a fusion bundle over $LM$, and $x\in M$. We denote by $\px Mx \subset PM$ the subspace of those paths that start at $x$. Then, there is a diffeological bundle gerbe $\un_x(P,\lambda)$ over $M$ consisting of
\begin{enumerate}[(i)]

\item 
the subduction $\ev_1: \px Mx \to M: \gamma \mapsto \gamma(1)$.

\item
the diffeological principal $S^1$-bundle $\ell^{*}P$ over $\px Mx^{[2]}$.

\item
the  gerbe product $\lambda$ on $\ell^{*}P$.

\end{enumerate}
This defines a \emph{regression functor}
\begin{equation*}
\un_x : \ufusbun {LM} \to \hc 1\ugrb M\text{.}
\end{equation*}
The main theorem of the transgression-regression machine is that regression is inverse to transgression, in the following sense: 

\begin{theorem}
\label{th:transreg}
Let $M$ be a connected smooth manifold. Then, the diagram
\begin{equation*}
\alxydim{@R=\diagR}{ & \ufusbun {LM} \ar[dr]^{\un_x} & \\ \hc 1 \ugrbcon M \ar[ur]^{L} \ar[rr] &&  \hc 1 \ugrb M }
\end{equation*}
of functors, which has on the bottom the functor that forgets connections and embeds bundle gerbes into diffeological bundle gerbes, is commutative up to a canonical natural equivalence. 
\end{theorem}

Theorem \ref{th:transreg} is proved in \cite[Section 6.1]{waldorf10} by constructing for each bundle gerbe $\mathcal{G}$ with connection over $M$  a 1-isomorphism
\begin{equation*}
\mathcal{A}_{\mathcal{G},y}: \mathcal{G} \to \un_x(L{\mathcal{G}},\lambda_{\mathcal{G}})\text{.}
\end{equation*}
This 1-isomorphism depends on the additional choice of a lift $y\in Y$ of the base point $x\in M$ along the surjective submersion of the bundle gerbe $\mathcal{G}$. Different choices of $y$ lead to 2-isomorphic 1-isomorphisms, $\mathcal{A}_{\mathcal{G},y} \cong \mathcal{A}_{\mathcal{G},y'}$. Under the operation $\hc 1$, these 2-isomorphisms become equalities; the resulting morphism $\hc 1 \mathcal{A}_{\mathcal{G},y}$ is thus independent of the choice of $y$. 

\begin{remark}
Transgression and regression can be made an equivalence of categories by either incorporating the connections on the side of the fusion bundles, or dropping the connections on the side of the gerbes; see the main theorems of \cite{waldorf10,waldorf11}. 
\end{remark}

Transgression and regression can be promoted to a multiplicative setting, i.e. with \emph{multiplicative} bundle gerbes (with connection) on the left hand side. On the loop space side we need:

\begin{definition}
\label{def:fusionextension}
A \emph{fusion extension} of $LG$ is a central extension
\begin{equation*}
1 \to S^1 \to P \to LG \to 1
\end{equation*}
of diffeological groups together with a multiplicative fusion product $\lambda$ on $P$.
\end{definition}

Here it is important that the evaluation map $\ev:PG^{[2]} \to G$, as well as path composition and inversion  are group homomorphisms.
In particular, the map $\ell: PG^{[2]} \to LG$ is a group homomorphism. The multiplicativity condition for the fusion product  is  that
\begin{equation*}
\lambda(q_{12} \otimes q_{23}) \cdot \lambda(q_{12}' \otimes q_{23}') = \lambda(q^{}_{12}q_{12}' \otimes q^{}_{23}q_{23}')
\end{equation*}
for all elements $q_{ij}\in P_{\ell(\gamma_i,\gamma_j)}$ and $q'_{ij}\in P_{\ell(\gamma_i',\gamma_j')}$ and all $(\gamma_1,\gamma_2,\gamma_3),(\gamma_1',\gamma_2',\gamma_3') \in PG^{[3]}$.

One can show that transgression sends a multiplicative bundle gerbe with connection to a fusion extension \cite[Section 1.3]{waldorf10}.
Here, it will be more important to look at \emph{regression}. With  the base point $1\in G$ understood,  a fusion bundle $(P,\lambda)$ over $LG$ regresses to a  diffeological bundle gerbe $\mathscr{R}(P,\lambda)$.  It is easy to check that the additional structure of a fusion extension (the group structure on $P$) makes $\un(P,\lambda)$   a strictly multiplicative, diffeological bundle gerbe.

\begin{remark}
Transgression and regression  can be seen as a functorial \emph{strictification}
\begin{equation*}
\bigset{8.3em}{Multiplicative bundle gerbes with connection over $G$}
\to 
\bigset{4.5em}{Fusion extensions of $LG$} 
\to
\bigset{9.8em}{Strictly multiplicative, diffeological bundle gerbes over $G$}\text{.}
\end{equation*} 
If $G$ is compact and simple, so that $\h^4(BG,\Z)\to\h^3(G,\Z)$ is injective, it follows from Theorem \ref{th:transreg} that this strictification preserves the characteristic class in $\h^4(BG,\Z)$. 
\end{remark}

\section{The Mickelsson Product}

\label{sec:mickelsson}

In this section we suppose that $G$ is compact, connected and simply-connected, for example $G\eq \spin n$ for $n > 2$. 
We consider the differential forms
\begin{equation}
\label{eq:forms}
H :=\frac{1}{6}\left \langle \theta \wedge [\theta \wedge \theta]  \right \rangle \in \Omega^3(G)
\quand
\rho := \frac{1}{2}\left \langle  \pr_1^{*}\theta \wedge \pr_2^{*}\bar\theta  \right \rangle \in \Omega^2(G \times G)\text{,}
\end{equation}
where $\theta$ and $\bar\theta$ are the left and right invariant Maurer-Cartan forms on $G$, respectively, and $\left \langle -,-  \right \rangle$ is an invariant bilinear form on the Lie algebra $\mathfrak{g}$ of $G$.  The  forms $H$ and $\rho$ satisfy the identities
\begin{equation}
\label{eq:formprop}
\mathrm{d}H=0 \quomma \Delta H = \mathrm{d}\rho \quand \Delta\rho = 0\text{,}
\end{equation}
where $\Delta: \Omega^q(G^k) \to \Omega^q(G^{k+1})$ is the alternating sum over the pullbacks along the face maps of the nerve of   $\B G$. Hence, the second and third equation become (in the notation of Section \ref{sec:multgrb}) 
\begin{equation*}
\pr_1^{*}H - m^{*}H + \pr_2^{*}H=\mathrm{d}\rho
\quand
\rho_{1,2} - \rho_{2,3} + \rho_{12,3} - \rho_{1,23} = 0\text{.}
\end{equation*}

Suppose  $\mathcal{G}$ is a bundle gerbe over $G$ with connection of curvature $H$. The \emph{Mickelsson product}
\begin{equation*}
\mp : L{\mathcal{G}} \times L{\mathcal{G}} \to L{\mathcal{G}}
\end{equation*}
on the transgression of $\mathcal{G}$ is defined as follows \cite{mickelsson1}; see \cite[Section 3.1]{waldorf5} and \cite[Theorem 6.4.1]{brylinski1}. First of all, we recall that the connection on the bundle gerbe $\mathcal{G}$ determines a \emph{surface holonomy} $\hol{\mathcal{G}}\varphi \in S^1$ for every closed oriented surface $\Sigma$ and a smooth map $\varphi:\Sigma \to G$. In its application to two-dimensional field theories, the surface holonomy provides the Feynman amplitude of the so-called Wess-Zumino term \cite{gawedzki3}. If the surface $\Sigma$ has a boundary one has to impose a boundary condition in order to keep the  holonomy well-defined. The boundary condition may be provided by a trivialization $\mathcal{T}$ of $\varphi^{*}\mathcal{G}|_{\partial\Sigma}$ \cite{carey2}; the surface holonomy in this case is denoted by $\mathcal{A}_{\mathcal{G}}(\varphi,\mathcal{T})$. We refer to \cite[Section 3.3]{waldorf10} for a detailed treatment with more references.

For loops $\tau,\tau'\in LG$, let $\mathcal{T},\mathcal{T}'$ be trivializations of $\tau^{*}\mathcal{G}$ and $\tau^{\prime *}\mathcal{G}$; these represent elements in $L\mathcal{G}$ over $\tau$ and $\tau'$, respectively. We choose extensions $\varphi,\varphi'\maps D^2 \to G$ of $\tau$ and $\tau'$ from the circle $S^1$ to its bounding disc $D^2$; these exist because $G$ is simply connected. The pointwise product $\tilde \varphi \df \varphi\varphi'$ is a similar extension of $\tilde\tau \df \tau\tau'$. We choose any trivialization $\widetilde{\mathcal{T}}$ of $\tilde\tau^{*}\mathcal{G}$. Finally, we consider  the combined map $\Phi \df (\varphi,\varphi')\maps D^2 \to G \times G$. Then, we define the Mickelsson product by
\begin{equation}
\label{def:prod}
\mathcal{T} \mp \mathcal{T}' := \widetilde{\mathcal{T}} \cdot \mathscr{A}_{\mathcal{G}}(\varphi,\mathcal{T})^{-1} \cdot \mathscr{A}_{\mathcal{G}}(\varphi',\mathcal{T}') ^{-1}\cdot \mathscr{A}_{\mathcal{G}}(\tilde\varphi,\widetilde{\mathcal{T}}) \cdot \exp \left ( \int_{D^2} \Phi^{*}\rho \right )\text{,}
\end{equation} 
where $\cdot$ denotes the action of $S^1$ on the element  $\widetilde{\mathcal{T}} \in L\mathcal{G}$. 

\begin{lemma}
Definition \erf{def:prod} is independent of all choices and turns $L\mathcal{G}$ into a central extension
\begin{equation*}
1 \to S^1 \to L\mathcal{G} \to LG \to 1
\end{equation*} 
of diffeological groups.
\end{lemma}

\begin{proof}
Suppose $\widetilde{\mathcal{T}}_1$ and $\widetilde{\mathcal{T}}_2$ are two choices of trivializations.  By \erf{eq:torsor} there exists a principal $S^1$-bundle $P$ with flat connection over $S^1$ such that $\widetilde{\mathcal{T}}_2 \cong \widetilde{\mathcal{T}}_1 \otimes P$. The associated surface holonomies satisfy  $\mathscr{A}_{\mathcal{G}}(\tilde\varphi,\widetilde{\mathcal{T}}_2) = \mathscr{A}_{\mathcal{G}}(\tilde\varphi,\widetilde{\mathcal{T}_1}) \otimes \mathrm{Hol}_P(S^1)^{-1}$  \cite[Lemma 3.3.2 (a)]{waldorf10}; this shows that \erf{def:prod} is independent of the choice of $\widetilde{\mathcal{T}}$. Suppose further that $(\varphi_1^{},\varphi_1')$ and $(\varphi_2^{},\varphi_2')$ are two choices of extensions of  $\tau$, $\tau'$. We consider the 2-sphere $S^2=D^2\#D^2$ as glued together from two discs, equipped with piecewise defined maps $\alpha := \varphi_1 \# \varphi_2$, $\alpha' := \varphi_1' \# \varphi_2'$ and $\tilde\alpha := \tilde \varphi_1 \# \tilde\varphi_2$, where $\tilde\varphi_k := \varphi_k^{} \varphi_k'$. The gluing law for surface holonomies  \cite[Lemma 3.3.2 (c)]{waldorf10} implies
\begin{equation}
\label{eq:gluing}
 \mathscr{A}_{\mathcal{G}}(\varphi_2,\mathcal{T}) = \mathscr{A}_{\mathcal{G}}(\varphi_1,\mathcal{T}) \cdot \mathrm{Hol}_{\mathcal{G}}(\alpha)\text{,}
\end{equation}
and analogous formulae with primes and tildes. Further, we consider the map $\Phi := \Phi_1 \# \Phi_2$ with $\Phi_i := (\varphi_i^{},\varphi_i')$. The identity $\Delta H = \mathrm{d}\rho$ implies \cite{gawedzki5,gawedzki9} the \emph{Polyakov-Wiegmann formula}
\begin{equation}
\label{eq:pw}
\mathrm{Hol}_{\mathcal{G}}(\alpha\alpha') = \mathrm{Hol}_{\mathcal{G}}(\alpha) \cdot \mathrm{Hol}_{\mathcal{G}}(\alpha') \cdot \exp \left ( \int_{S^2} \Phi^{*}\rho \right )\text{.}
\end{equation}
Formulas \erf{eq:gluing} and \erf{eq:pw} prove that  \erf{def:prod}   is independent of the choice of the extensions $\varphi$ and $\varphi'$.

Associativity of $\mp$ follows from $\Delta\rho=0$; smoothness  from the smoothness of the surface holonomy $\mathscr{A}_{\mathcal{G}}$  \cite[Lemma 4.2.2]{waldorf10}. The construction of a unit and of inverses is straightforward. 
Thus, $L\mathcal{G}$ is a diffeological group and also a principal $S^1$-bundle over $LG$, i.e  a central extension.
\end{proof}

Next we recall from Section \ref{sec:transgressionregression} that $L\mathcal{G}$ carries a fusion product $\lambda_{\mathcal{G}}$.

\begin{lemma}
\label{lem:canmult}
The fusion product $\lambda_{\mathcal{G}}$ is multiplicative with respect to the Mickelsson product. 
\end{lemma}

\begin{proof}
First we  mention the following general fact, for a bundle gerbe $\mathcal{G}$ with connection over a compact, simply-connected manifold $M$. Suppose $(\gamma_1,\gamma_2,\gamma_3) \in PM^{[3]}$. Since $M$ is simply-connected, there exists a smooth path $\Gamma: [0,1] \to PM^{[3]}$ such that $\Gamma(0)\eq (\gamma_1,\gamma_2,\gamma_3)$, and $\Gamma(1)$ is a triple of identity paths at some point in $M$. The paths $\varphi_{ij} \df \ell \circ \pr_{ij} \circ \Gamma$ in $LM$ can be regarded as extensions of the loops $\tau_{ij}:=\ell(\gamma_i,\gamma_j)$ to the disc. Then, \cite[Proposition 4.3.4]{waldorf10} implies that
\begin{equation}
\label{eq:surfmult}
\mathscr{A}_{\mathcal{G}}(\varphi_{12},\mathcal{T}_{12}) \cdot \mathscr{A}_{\mathcal{G}}(\varphi_{23},\mathcal{T}_{23}) = \mathscr{A}_{\mathcal{G}}(\varphi_{13},\mathcal{T}_{13})
\end{equation}
for any triple of trivializations $\mathcal{T}_{ij}$ of $\tau_{ij}^{*}\mathcal{G}$ satisfying $\lambda_{\mathcal{G}}(\mathcal{T}_{12}\otimes\mathcal{T}_{23})=\mathcal{T}_{13}$.

Now suppose $(\gamma_1,\gamma_2,\gamma_3),(\gamma_1',\gamma_2',\gamma_3') \in PG^{[3]}$ and $\mathcal{T}_{ij},\mathcal{T}_{ij}'$ are trivializations over $\tau_{ij},\tau_{ij}'$ such that  $\lambda_{\mathcal{G}}(\mathcal{T}_{12}\otimes\mathcal{T}_{23})=\mathcal{T}_{13}$ and  $\lambda_{\mathcal{G}}(\mathcal{T}_{12}'\otimes\mathcal{T}_{23}')=\mathcal{T}_{13}'$. We choose paths $\Gamma$, $\Gamma'$ as above, and  extract the extensions $\varphi_{ij}$, $\varphi_{ij}'$ each satisfying \erf{eq:surfmult}. The product $\widetilde{\Gamma} := \Gamma \cdot \Gamma'$ produces the extensions $\tilde\varphi_{ij} = \varphi_{ij}\varphi_{ij}'$ also satisfying \erf{eq:surfmult}. For the combined maps $\Phi_{ij}=(\varphi_{ij},\varphi_{ij}')$ we have by construction
\begin{equation}
\label{eq:integrals}
\int_{D^2} \Phi^{*}_{13}\rho = \int_{D^2} \Phi^{*}_{12}\rho + \Phi_{23}^{*}\rho\text{.}
\end{equation}
Define $\widetilde{\mathcal{T}_{12}} := \mathcal{T}_{12} \mp \mathcal{T}_{12}'$ and $\widetilde{\mathcal{T}_{23}} := \mathcal{T}_{23} \mp \mathcal{T}_{23}'$, i.e. these are trivializations that satisfy via \erf{def:prod}
\begin{equation}
\label{eq:mickelssontriv}
\mathscr{A}_{\mathcal{G}}(\varphi_{ij},\mathcal{T}_{ij}) \cdot \mathscr{A}_{\mathcal{G}}(\varphi'_{ij},\mathcal{T}_{ij}') = \mathscr{A}_{\mathcal{G}}(\tilde\varphi_{ij},\widetilde{\mathcal{T}}_{ij}) \cdot \exp \left ( \int_{D^2} \Phi_{ij}^{*}\rho \right )\text{.}
\end{equation}
The multiplicativity we have to show is now equivalent to the identity
\begin{equation*}
\mathcal{T}_{13} \mp \mathcal{T}_{13}' = \lambda_{\mathcal{G}}(\widetilde{\mathcal{T}_{12}}\otimes  \widetilde{\mathcal{T}_{23}})\text{.}
\end{equation*}
It follows  from \erf{eq:surfmult}, \erf{eq:integrals} and \erf{eq:mickelssontriv} upon computing the left hand side with  $\widetilde{\mathcal{T}_{13}} := \lambda_{\mathcal{G}}(\widetilde{\mathcal{T}_{12}} \otimes \widetilde{\mathcal{T}_{23}})$.
\end{proof}

Summarizing, we obtain:

\begin{theorem}
\label{th:fusext}
Let $G$ be a compact, connected, simply-connected Lie group, and let $\mathcal{G}$ be a bundle gerbe over $G$ with connection of curvature $H$. Then,
the Mickelsson product equips the transgression $L\mathcal{G}$ with the structure of a fusion extension of $LG$.
\end{theorem}

\section{The Construction of String 2-Group Models}

\label{sec:string}

In this section we  consider a  compact,  simple, simply-connected Lie group $G$ such as $\spin n$ for $n=3$ or $n> 4$. We  briefly review the  \quot{basic} bundle gerbe $\gbas$ over $G$ whose Dixmier-Douady class generates $\h^3(G,\Z)\cong \Z$, following Gaw\c edzki-Reis \cite{gawedzki1,gawedzki2}, Meinrenken \cite{meinrenken1}, and Nikolaus \cite{nikolaus1}.   
  
We choose a Weyl alcove $\mathfrak{A}$ in the dual $\mathfrak{t}^{*}$ of the Lie algebra of a maximal torus of $G$. For these exist canonical choices  \cite[Section 4]{gawedzki2}. The alcove $\mathfrak{A}$ parameterizes conjugacy classes of $G$ in terms of a continuous map $q: G \to \mathfrak{A}$. We denote by $\mathfrak{A}_{\mu} := \mathfrak{A}\setminus f_{\mu}$ the complement of the closed face $f_{\mu}$ opposite to a vertex $\mu$ of $\mathfrak{A}$. 
The preimages $U_\mu$ of $\mathfrak{A}_{\mu}$ under $q$ form a cover of $G$ by open sets. We denote by $G_\mu$ the centralizer of $\mu$ in $G$ under the coadjoint action. These centralizer groups come with central $S^1$-extensions $\hat G_{\mu}$ which are trivial if and only if $G_{\mu}$ is simply-connected. Each open set $U_{\mu}$ supports a smooth map $\rho_\mu: U_{\mu} \to G/G_{\mu}$, and thus the principal $G_{\mu}$-bundle $P_{\mu} := \rho_i^{*}G$. The problem of lifting the structure group of $P_{\mu}$ from $G_{\mu}$ to $\hat G_{\mu}$ defines a \emph{lifting bundle gerbe} $\mathcal{L}_{\mu}$ over $U_{\mu}$. These local lifting bundle gerbes glue together and yield the basic gerbe $\gbas$. Further, each  $\mathcal{L}_{\mu}$ can be equipped with a connection, and the glued connection on $\gbas$ has curvature $H$, for a certain normalization of the bilinear form $\left \langle -,-  \right \rangle$ in \erf{eq:forms}.

The transgression  $L\gbas$ is a fusion extension of $LG$ (Theorem \ref{th:fusext}), so that the multiplicative regression functor  of Section \ref{sec:transgressionregression} produces a  strictly multiplicative, diffeological bundle gerbe 
\begin{equation*}
\mathcal{R} :=\un(L\gbas,\lambda_{\gbas})
\end{equation*}
over $G$. 
We may now optionally proceed in the following two ways: 

1.) 
Theorem \ref{th:transreg} shows that $\mathcal{R} \cong \gbas$; whence the  class of $\mathcal{R}$ generates $\h^3(G,\Z)\cong \h^4(BG,\Z)$. Thus, the 2-functor \erf{eq:cextstrictdiff} produces a central, strict, diffeological 2-group extension
\begin{equation*}
\B S^1 \to \grpd{\mathcal{R}} \to \idmorph{G}
\end{equation*}
with the same class, so that, for $G=\spin n$, $\grpd{\mathcal{R}}$ is a 2-group model for $\str n$. Let us summarize the structure of $\grpd{\mathcal{R}}$ by assembling the various constructions: its space of objects is $P_1G$ and its space of morphisms is $\ell^{*}L\gbas = P_1G^{[2]} \lli{\ell}\times_{\pr} L\gbas$,  composition is the fusion product $\lambda_{\gbas}$, and  multiplication is the Mickelsson product. We remark that $\grpd{\mathcal{R}}$ has  (essentially) the same objects and morphisms as the model of \cite{baez9}, but  the composition is defined in \cite{baez9} using the multiplication (the Mickelsson product) and here using the fusion product.

2.) 
Theorem \ref{th:transreg} not only shows that $\mathcal{R} \cong \gbas$, it also provides a distinguished 1-isomorphism
\begin{equation*}
\mathcal{A}_{\gbas,y}: \gbas \to \mathcal{R}\text{,}
\end{equation*}
where $y\in Y$ is an element in the surjective submersion of $\gbas$ that projects to $1\in G$. In the  con\-struc\-tion of $\gbas$ outlined above there is a such an element: the identity element $1 \in G$ lies in the open set $U_0$ associated to the origin $0 \in \mathfrak{g}^{*}$. Accordingly, its stabilizer is $G_0 = G$, and $P_0$ is the trivial principal $G$-bundle over $U_0$. As such, it has a canonical element $p = (1,1) \in P_0 = U_0 \times G$. In the gluing construction of the local lifting  gerbes $\mathcal{L}_{\mu}$ the surjective submersion $\pi:Y \to G$ of $\gbas$ is the disjoint union of total spaces $P_{\mu}$ of the submersions of $\mathcal{L}_{\mu}$; thus, $p\in Y$.
Now, the multiplicative structure on $\mathcal{R}$ can be \quot{pulled back} to $\gbas$ along $\mathcal{A}_{\mathcal{G},p}$.

The result is a \emph{diffeological} multiplicative structure on the \emph{finite-dimensional} bundle gerbe $\gbas$. Its 1-isomorphism $\mathcal{M}$ involves a certain subduction $\chi: Z \to Y' :=Y_{1,2} \times_{G \times G} Y_{12}$, where $Y_{1,2}$ and $Y_{12}$ are the smooth manifolds we have encountered in Section \ref{sec:multgrb}. It further involves a diffeological principal $S^1$-bundle $Q$ over $Z$.  General bundle gerbe theory \cite[Theorem 1]{waldorf1} shows that $Q$ descends along $\chi:Z \to Y'$. But a diffeological principal $S^1$-bundle over a smooth manifold is automatically smooth \cite[Theorem 3.1.7]{waldorf9}. This defines a new, smooth 1-isomorphism $\mathcal{M}'$. Both steps are functorial so that the associator $\alpha$ for $\mathcal{M}$ descends to an associator $\alpha'$ for $\mathcal{M}'$. Since smooth manifolds embed fully and faithfully into diffeological spaces, it follows that  $\alpha'$ is smooth. Thus, $(\gbas,\mathcal{M}',\alpha')$ is a smooth, multiplicative bundle gerbe over $G$ whose class generates $\h^4(BG,\Z)$. Under the 2-functor  \erf{eq:diagram} it hence yields a smooth, finite-dimensional Lie 2-group extension $\grpd{\gbas}$ of $G$ by $\B S^1$ of the same class. In particular, for $G=\spin n$,  it is a  2-group model for $\str n$.

\kobib{../../bibliothek/tex}

\end{document}